\documentclass{amsart}
\usepackage{amssymb}
\usepackage{amsthm,amsmath}

\usepackage{graphicx,color}
\usepackage[all]{xy}

\usepackage{enumerate}
\usepackage{ifpdf}
\usepackage{mathtools}

\newtheorem{theorem}{Theorem}[section]

\newtheorem{lemma}[theorem]{Lemma}
\newtheorem{proposition}[theorem]{Proposition}
\theoremstyle{definition}
\newtheorem{definition}[theorem]{Definition}
\theoremstyle{remark}

\numberwithin{equation}{section}
\allowdisplaybreaks

\newcommand{\M}{\mathcal{M}}
\newcommand{\B}{\mathcal{B}}
\newcommand{\Z}{\mathbb{Z}}
\newcommand{\ra}{\rightarrow}
\renewcommand{\bar}{\overline}
\renewcommand{\loop}{\Omega}
\newcommand{\cat}{\mathcal}
\newcommand{\id}{1}
\newcommand{\aut}{\operatorname{Aut}}

\begin{document}




\title{The braidings in the mapping class groups of surfaces}


\author{Yongjin Song}
\address{Departments of Mathematics \\ Inha University \\ Incheon 402-751, Korea}
\email{yjsong@inha.ac.kr}

\begin{abstract}
The disjoint union of mapping class groups of surfaces forms a braided monoidal category $\mathcal M$, as the disjoint union of the braid groups $\mathcal B$ does. We give a concrete, and geometric meaning of the braiding $\beta_{r,s}$ in $\M$. Moreover, we find a set of elements in the mapping class groups which correspond to the standard generators of the braid groups. Using this, we obtain an obvious map $\phi:B_g\ra\Gamma_{g,1}$. We show that this map $\phi$ is injective and nongeometric in the sense of Wajnryb. Since this map extends to a braided monoidal functor $\Phi : \mathcal B \rightarrow \mathcal M$, the integral homology homomorphism induced by $\phi$ is trivial in the stable range.
\end{abstract}



\maketitle

\section{Introduction}

Let $\Gamma_{g,1}$ be the mapping class group of the surface $S_{g,1}$, the compact orientable surface of genus $g$ with one boundary component obtained by deleting a disk from the closed surface $S_g$. Let $\M=\coprod_{g\geq 0}\Gamma_{g,1}$ be the disjoint union of $\Gamma_{g,1}$'s, then it is regarded as the category whose objects are nonnegative integers. Here, $\hom_{\M}(g,g)=\Gamma_{g,1}$ and there is no morphism between distinct integers. The monoid structure, called F-product in $\M$, is induced by the pair-of-pants connected sum. The braid structure of $\M$ was given in \cite{Song} in terms of the actions of braidings on the fundamental group of $S_{g,1}$. Recall that $\Gamma_{g,1}$ may be identified with the subgroup of the automorphism group of $\pi_1 S_{g,1}$ that consists of the automorphisms fixing the fundamental relator $R=[y_1,x_1]\cdots [y_g,x_g]$. Here $\pi_1 S_{g,1}$ is a free group on $2g$ generators $x_1,y_1,x_2,y_2,\ldots,x_g,y_g$ and the fundamental relator represents a loop along the boundary of $S_{g,1}$. The F-product on $\M$ may be identified with the operation taking the free product of automorphisms.

The action of the braiding $\beta_{r,s}:r\otimes s\ra s\otimes r$, which is an element of $\Gamma_{r+s,1}$, on the free group on $x_1,y_1,x_2,y_2,\ldots,x_{r+s},y_{r+s}$ was given in \cite{Song}. The braiding could be also expressed as a product of the standard Dehn twists, the generators of $\Gamma_{g,1}$.

The disjoint union of braid groups $\mathcal B=\coprod_{g\geq 0}B_g$ forms a braided monoidal category. Let $\sigma_1,\ldots,\sigma_{g-1}$ be the standard generators of $B_g$, where $\sigma_i$ crosses the $i$-th and $i+1$-st strings. The monoid structure is induced by the juxtaposition and the braiding $(\sigma_{r,s}:r\otimes s\ra s\otimes r)\in B_{r+s}$ is the crossing of the front $r$ strings and the rear $s$ strings. We can express $\sigma_{r,s}$ in terms of $\sigma_i$'s as follows:
\[\sigma_{r,s}=(\sigma_r\sigma_{r-1}\cdots\sigma_1)(\sigma_{r+1}\sigma_r\cdots\sigma_2)\cdots(\sigma_{r+s-1}\sigma_{r+s-2}\cdots\sigma_s).\]

In this paper, we give a concrete, and geometric meaning of the braiding $\beta_{r,s}$ in the braided monoidal category $\M$. Moreover, we find a set of elements $\beta_1,\ldots,\beta_{g-1}$ of $\Gamma_{g,1}$ which are analogous to  the standard generators of $B_g$. Using these generators we can define an obvious map $\phi:B_g\ra\Gamma_{g,1}, \sigma_i\mapsto\beta_i$. An embedding of a braid group into a mapping class group is said to be {\it geometric} by Wajnryb (\cite{W2}, \cite{W3} ) if it maps the standard generators of the braid group to Dehn twists. We show that the map $\phi:B_g\rightarrow\Gamma_{g,1}$ is injective and is, moreover, {\it nongeometric}.

The natural embedding $\phi:B_g\ra\Gamma_{g,1}$ can be extended to a braided monoidal functor $\Phi:\mathcal B\ra\M$. Namely, $\Phi$ preserves the braidings. As a consequence we show that $\phi:B_g\ra\Gamma_{g,1}$ induces a map of double loop spaces $B\phi^+:BB_\infty^+\ra B\Gamma_\infty^+$ and hence the homology homomorphism $\phi_*:H_*(B_\infty;\Z)\ra H_*(\Gamma_\infty;\Z)$ induced by $\phi$is trivial, where $B_\infty=\varinjlim B_g$ and $\Gamma_\infty=\varinjlim\Gamma_{g,1}$.

The classifying space of a braided monoidal category naturally gives rise to a double loop space(\cite{Fie}). More precisely speaking, the group completion of the classifying space of a braided monoidal category is homotopy equivalent to a double loop space. For $\mathcal B=\coprod_{g\geq 0} B_g$, we have
\[\bar{B\mathcal B}=\loop B(\coprod_{g\geq 0}BB_g)\simeq \Z_+\times BB_\infty^+ \simeq \loop^2S^2,\]
where $\bar{B\mathcal B}$ means the group completion of the classifying space of $\mathcal B$, and $+$ means the Quillen's plus construction.

The braided monoidal functor $\Phi:\mathcal B\ra\M$ induces a map of double loop spaces $\bar{B\mathcal B}\ra\bar{B\M}$. In other words, $B\phi^+:BB_\infty^+\ra B\Gamma_\infty^+$ is a map of double loop spaces. Since every map $\Psi:BB_\infty^+\ra B\Gamma_\infty^+$ of double loop spaces is null-homotopic (\cite{S-T}, Lemma 5.3), so is $B\phi^+$. This proves that the homomorphism $\phi_*:H_*(B_\infty;\Z)\ra H_*(\Gamma_\infty;\Z)$ is trivial, since the plus construction does not change the homology.

When we first think about the geometric meaning of he braiding $\beta_{r,s}$, since it interchanges the front $r$ genus holes and rear $s$ genus holes of $S_{r+s,1}$, we think that it would be a $180^\circ$ rotation around the boundary  of the surface. Here we consider the boundary as the waist of the surface with two arms; the left arm is the front $r$ genus holes and the right arm is the rear $s$ genus holes. However, just a $180^\circ$ rotation around the waist is not sufficient. We also need half-twists of two arms. Here, it is important to determine the directions of these half-twists. There are four choices of directions; two choices on the waist and two choices on the arms. By analyzing the simplest case, $\beta_{1,1}$ for $S_{2,1}$, we show that the braiding $\beta_{r,s}$ is the composite of two half-twists; the {\it reverse} half Dehn twist along the boundary and half Dehn twists along the shoulders of two arms. The key part of this paper is to figure out a correct geometric meaning of $\beta_{1,1}$ as a self-homeomorphism of $S_{2,1}$ and extend this to the general case, $\beta_{r,s}$ in $\Gamma_{r+s,1}$.

\section{Braided monoidal category}

\begin{figure}[hbt]
\begin{center}
  \includegraphics{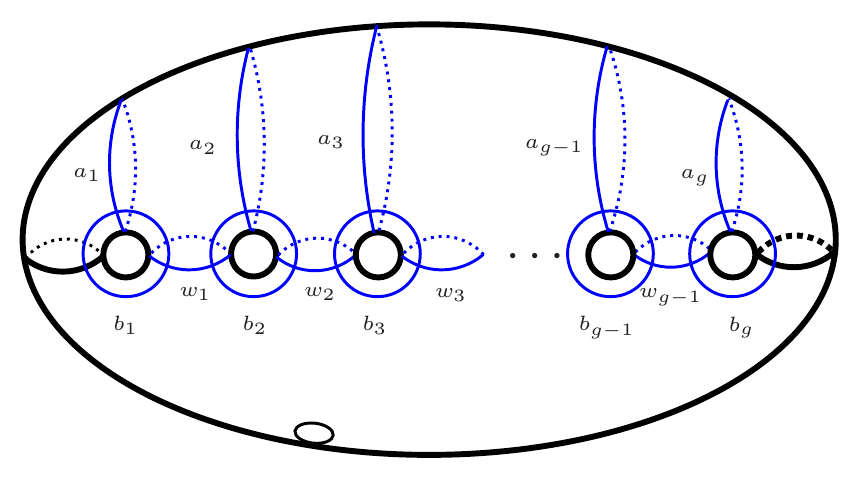}\\
  \caption{Standard Dehn twists of $\Gamma_{g,1}$}\label{fig1}
\end{center}
\end{figure}

Let $S_{g,1}$ be a compact connected orientable surface of genus $g$ with one boundary component. The mapping class group $\Gamma_{g,1}$ is the group of isotopy classes of orientation preserving self-homeomorphisms of $S_{g,1}$ fixing the boundary of $S_{g,1}$ pointwise. $\Gamma_{g,1}$ is generated by the standard Dehn twists $a_1,\ldots,a_g,b_1,\ldots,b_g,w_1,\ldots,w_{g-1}$ (Figure~\ref{fig1}).

\begin{definition}\label{def1}
  A {\it braided monoidal category} is a monoidal category $(\cat C,\otimes, I)$ together with a family of natural commutativity isomorphisms
  \[\beta_{A,B}:A\otimes B\ra B\otimes A\]
  called {\it braidings}, natural in both variables, such that for all objects $A,B,C$ in $\cat C$,
  \begin{enumerate}[(a)]
  \item $\beta_{A,I}=\beta_{I,A}=\id_A$
  \item $\beta_{A\otimes B,C}=(\beta_{A,C}\otimes\id_B)\circ(\id_A\circ\beta_{B,C})$
  \item $\beta_{A,B\otimes C}=(\id_B\otimes\beta_{A,C})\circ(\beta_{A,B}\otimes\id_C)$.
  \end{enumerate}
\end{definition}
The equations in $(b),(c)$ imply the Yang-Baxter equation:
\[(\id_C\otimes\beta_{A,B})(\beta_{A,C}\otimes\id_B)(\id_A\otimes\beta_{B,C})=(\beta_{B,C}\otimes\id_A)(\id_B\otimes\beta_{A,C})(\beta_{A,C}\otimes\id_C).\]

The disjoint union $\M=\coprod_{g\geq 0}\Gamma_{g,1}$ forms a category whose objects are nonnegative integers and morphisms satisfy
\[\hom(g,h)=\left\{
  \begin{array}{cl}
    \Gamma_{g,1} & \mathrm{\ if\ }g=h, \\ \phi & \mathrm{\ if\ }g\neq h.
  \end{array}
\right.\]
$\M$ has a monoid structure: the F-product
\[\Gamma_{g,1}\times \Gamma_{h,1}\ra\Gamma_{g+h,1}\]
is induced by extending two self-homeomorphisms on $S_{g,1}$ and $S_{h,1}$ to the surface $S_{g+h,1}$ obtained attaching a pair of pants (a sphere with three boundary components) to the surfaces $S_{g,1}$ and $S_{h,1}$ along the fixed boundary circles. We extend the identity map on the boundary to the whole pants.

It was shown in \cite{Song} that $\M$ is a braided monoidal category whose braidings are as follows: Let $x_1,y_1,x_2,y_2,\ldots,x_g,y_g$ be the generators of $\pi_1S_{g,1}$ which are represented by loops parallel to the Dehn twist loops $a_1,b_1,a_2,b_2,\ldots,a_g,b_g$, respectively, in Figure~\ref{fig1}. $\Gamma_{g,1}$ is identified with the subgroup of the automorphism group of $\pi_1S_{g,1}$ fixing the fundamental relator $R=[y_1,x_1][y_2,x_2]\cdots[y_g,x_g]$. The $(r,s)$-braiding $\beta_{r,s}:r\otimes s\ra s\otimes r$, which is an element of $\Gamma_{r+s,1}$, acts on the free group on $\{x_1,y_1,x_2,y_2,\ldots,x_g,y_g\}$ as follows:
\[\begin{aligned}
  x_1\mapsto R_sx_{s+1}R_s^{-1} & ,\ y_1\mapsto R_sy_{s+1}R_s^{-1} \\
  & \vdots\\
  x_r\mapsto R_sx_{s+r}R_s^{-1} & ,\ y_1\mapsto R_sy_{s+r}R_s^{-1} \\
  x_{r+1}\mapsto x_1 & ,\ y_{r+1}\mapsto y_1 \\
  x_{r+2}\mapsto x_2 & ,\ y_{r+1}\mapsto y_2 \\
  & \vdots\\
  x_{r+s}\mapsto x_s & ,\ y_{r+s}\mapsto y_s \\
\end{aligned}\]
where $R_s=[y_1,x_1]\cdots[y_s,x_s]$. We can easily check that the $(r,s)$-braidings fix the fundamental relator $R$ and satisfy the braid equations (b), (c) of Definition~\ref{def1}.

\begin{figure}[hbt]
\begin{center}
  \includegraphics{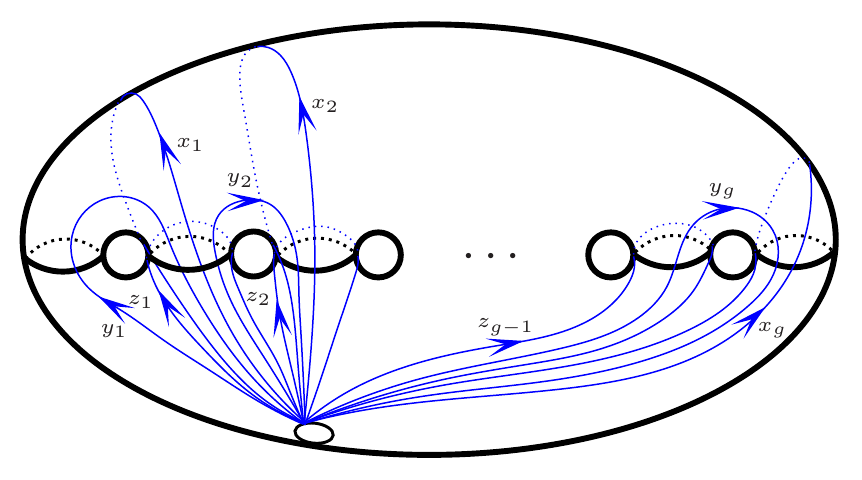}\\
  \caption{The generators $x_1,y_1,\ldots,x_g,y_g$ of $\pi_1S_{g,1}$}\label{fig2}
\end{center}
\end{figure}

Since the group completion of the classifying space of a braided monoidal category is  homotopy equivalent to a double loop space, we have
\begin{theorem}
  $\bar{B\M}=\loop B(\coprod_{g\geq 0}B\Gamma_{g,1})$ is homotopy equivalent to a double loop space. Here, $\bar{B\M}$ denotes the group completion of $B\M$.
\end{theorem}

\section{The geometric analysis of braidings}

In this section we give geometric meanings of braidings in the mapping class group. By using this result we construct a natural embedding $\phi:B_g\hookrightarrow\Gamma_{g,1}$, which is nongeometric. We show that this embedding induces the trivial homology homomorphism in the stable range in the integral coefficient.

The braiding $\beta_{r,s}\in\Gamma_{r+s,1}$ can be expressed as a product of the standard Dehn twists. First, the braiding $\beta_{1,1}$ was expressed in \cite{Song} in terms of the Dehn twists $a_1,a_2,b_1,b_2,w_1$:
\begin{lemma}\label{lemma1}
  The $(1,1)$-braiding $\beta_{1,1}$ in $\Gamma_{2,1}$ equals
  \begin{equation}
    \beta_{1,1}=(a_1b_1a_1)^4(a_2b_2(a_1b_1a_1)^{-1}w_1a_1b_1a_1^2b_1)^{-3}.\label{eq:one-one-brading}
  \end{equation}
\end{lemma}

The standard Dehn twists $a_1,\ldots,a_g,b_1,\ldots,b_g,w_1,\ldots,w_{g-1}$ act on $\pi_1S_{g,1}$, which is the free group generated by $\{x_1,y_1,x_2,y_2,\ldots,x_g,y_g\}$, as follows:
\[
\setlength{\arrayrulewidth}{0pt}
\begin{aligned}
  a_i:&y_i\mapsto y_ix_i^{-1}\\
  b_i:&x_i\mapsto x_iy_i\\
  w_i:&x_i\mapsto z_i^{-1}y_{i+1}x_{i+1}y_{i+1}^{-1}\\
  &y_i\mapsto y_iz_i\\
  &y_{i+1}\mapsto z_i^{-1}y_{i+1}
\end{aligned}
\]
where $z_i=x_i^{-1}y_{i+1}x_{i+1}y_{i+1}^{-1}$ and these automorphisms fix the generators that do not appear in the list. By using these actions we can check that the equation given in Lemma~\ref{lemma1} is correct, that is, $\beta_{1,1}$ acts on the free group on $\{x_1,y_1,x_2,y_2\}$ as
\[
\begin{aligned}
  x_1\mapsto[y_1,x_1]x_2[x_1,y_1]&,\ y_1\mapsto[y_1,x_1]y_2[x_1,y_1]\\
  x_2\mapsto x_1&,\ y_2\mapsto y_1
\end{aligned}
\]
as desired. Here the Dehn twists act on the right on $\pi_1S_{g,1}$.

Now we construct the $(r,s)$-braiding $\beta_{r,s}$ by extending $\beta_{1,1}$. Let
\[\beta_i=(a_ib_ia_i)^4(a_{i+1}b_{i+1}(a_ib_ia_i)^{-1}w_ia_ib_ia_i^2b_i)^{-3}\]
for $i=1,\ldots,g-1$, then these satisfy the braid relations. The $(r,s)$-braiding $\beta_{r,s}$ in the group $\Gamma_{g,1}$ defined by
\[\beta_{r,s}=(\beta_r\beta_{r-1}\cdots\beta_1)(\beta_{r+1}\beta_r\cdots\beta_2)\cdots(\beta_{r+s-1}\beta_{r+s-2}\cdots\beta_s)\]
is analogous to the $(r,s)$-braiding in the braid group $B_{r+s}$. Here we may think that the composition of two braids is a stacking (from up to) down as morphisms map from top to bottom. Through simple calculations we can check that this $(r,s)$-braiding $\beta_{r,s}$ acts on $\pi_1 S_{g,1}$ as described in Section 2.

We now investigate the geometric meaning of $\beta_{r,s}$. We first figure out what $\beta_{1,1}$ really is as a self-homeomorphism of $S_{2,1}$.

\begin{figure}[hbt]
\begin{center}
  \includegraphics{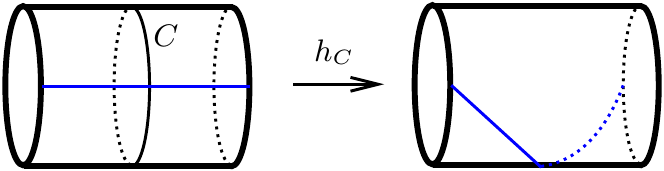}\\
  \caption{The half Dehn twist along $C$}\label{fig3}
\end{center}
\end{figure}

For a simple closed curve $C$ on $S_{g,1}$, the half Dehn twist $h_C$ is a self-homeomorphism of $S_{g,1}$ which is $180^\circ$ twist along a tubular neighborhood of $C$ in the same direction as the usual Dehn twist, as given in Figure~\ref{fig3}. Let $h'_C$ denote the reverse half Dehn twist along $C$, namely $h'_C=h_C^{-1}$ in $\Gamma_{g,1}$.

In this paper we sometimes abuse notations. For a simple closed curve $C$ on the surface, it would be just a curve for a (half) Dehn twist, but sometimes it stands for an element of the fundamental group of the surface. For example, Let $R$ be the closed curve along the boundary of $S_{g,1}$, then it could also mean the element $R=[y_1,x_1][y_2,x_2]\cdots[y_g,x_g]\in\pi_1S_{g,1}$.

\begin{figure}[hbt]
\begin{center}
  \includegraphics{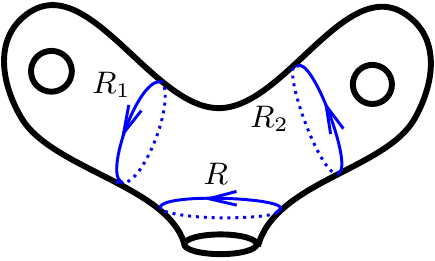}\\
  \caption{}\label{fig4}
\end{center}
\end{figure}

Let $R_1, R_2$ be closed curves on $S_{2,1}$ as given in Figure~\ref{fig4}. The arrows on the curves mean the directions of them as elements of $\pi_1S_{2,1}$.

\begin{figure}[hbt]
\begin{center}
  \includegraphics{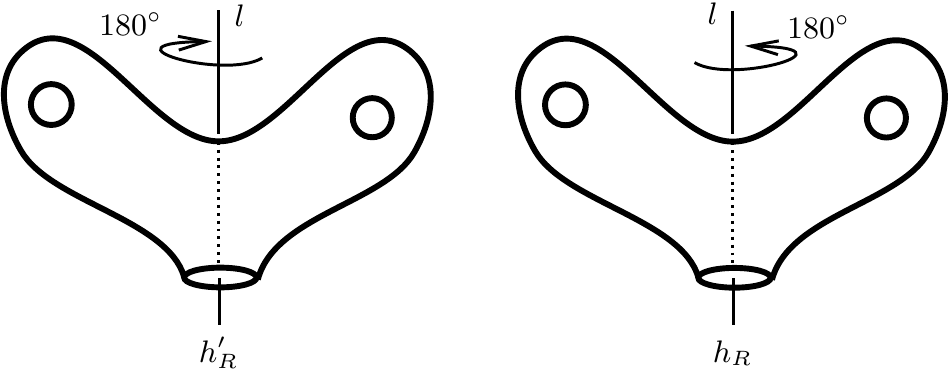}\\
  \caption{Two half Dehn twists along the boundary}\label{fig5}
\end{center}
\end{figure}

For the closed curve $R$ on $S_{2,1}$, the half Dehn twist $h_R$ and the reverse half Dehn twist $h_R'$ are illustrated in Figure~\ref{fig5}. We may think of them as $180^\circ$ rotations of the whole surface around the axis $l$.

We may regard $h_R$ and $h_R'$ as elements of the automorphism group of $\pi_1S_{2,1}=F_{\{x_1,y_1,x_2,y_2\}}$. The geometric calculations in this section not only play key roles in analyzing the geometry of braidings in $\Gamma_{g,1}$, but also are of interest in their own right.

\begin{proposition}
  The half Dehn twist $h_R$ acts on $\pi_1S_{2,1}$ as follows:
  \[\begin{aligned}
    x_1\mapsto R_2^{-1}x_2^{-1}&,\ y_1\mapsto y_2^{-1}R_2\\
    x_2\mapsto R^{-1}x_1^{-1}R_2&,\ y_2\mapsto R_2^{-1}y_1R.
  \end{aligned}\]
\end{proposition}

Note that $h_R$ maps $R_1$ to $R_2$, and $R_2$ to $R_2^{-1}R$.

\begin{proposition}
  The reverse half Dehn twist $h_R'$ acts on $\pi_1S_{2,1}$ as follows:
  \[\begin{aligned}
    x_1&\mapsto R_1x_2^{-1}R^{-1}=(R_1y_2)x_2^{-1}(R_1y_2)^{-1}\\
    y_1&\mapsto Ry_2^{-1}R_1^{-1}=(R_1y_2x_2)y_2^{-1}(R_1y_2x_2)^{-1}\\
    x_2&\mapsto y_1x_1^{-1}y_1^{-1}\\
    y_2&\mapsto R_1y_1^{-1}=(y_1x_1)y_1^{-1}(y_1x_1)^{-1}
  \end{aligned}\]
\end{proposition}

Note that $h_R'$ maps $R_1$ to $RR_1^{-1}$, and $R_2$ to $R_1$.

These lemmas are obtained through series of geometric calculations of actions of half Dehn twists on $\pi_1S_{2,1}=F_{\{x_1,y_1,x_2,y_2\}}$.

\begin{proposition}
  For the half Dehn twists $h_{R_1}$ and $h_{R_2}$, let $h_A=h_{R_1}\circ h_{R_2}$, called the half Dehn twists on two arms. Then $h_A$ acts on $\pi_1S_{g,1}$ as follows:
  \[
  \begin{aligned}
    x_1&\mapsto R_1^{-1}x_1^{-1}=(x_1y_1)x_1^{-1}(x_1y_1)^{-1}\\
    y_1&\mapsto y_1^{-1}R_1=x_1y_1^{-1}x_1^{-1}\\
    x_2&\mapsto R_2^{-1} x_2^{-1}=(x_2y_2)x_2^{-1}(x_2y_2)^{-1}\\
    y_2&\mapsto y_2^{-1}R_2=x_2y_2^{-1}x_2^{-1}
  \end{aligned}
  \]
\end{proposition}

\begin{proposition}
  For the reverse half Dehn twists on two arms $h_A'=h_{R_1}'\circ h_{R_2}'$ acts on $\pi_1S_{g,1}$ as follows:
  \[
  \begin{aligned}
    x_1&\mapsto x_1^{-1}R_1^{-1}=y_1x_1^{-1}y_1^{-1}\\
    y_1&\mapsto R_1y_1^{-1}=(y_1x_1)y_1^{-1}(y_1x_1)^{-1}\\
    x_2&\mapsto x_2^{-1}R_2^{-1}=y_2x_2^{-1}y_2^{-1}\\
    y_2&\mapsto R_2y_2^{-1}=(y_2x_2)y_2^{-1}(y_2x_2)^{-1}
  \end{aligned}
  \]
\end{proposition}

In the combination of the two half Dehn twists, the half Dehn twist along the boundary and the half Dehn twists on two arms, there are four choices of directions. For the braiding $\beta_{1,1}$, the right choice of directions is $h'_R$ and $h_A$.

\begin{theorem}
  The $(1,1)$-braiding $\beta_{1,1}\in\Gamma_{2,1}$ is the product (composite) of $h_R'$ and $h_A$ in $\Gamma_{2,1}$, that is, $\beta_{1,1}=h_R'\circ h_A$.
  \begin{proof}
    Since three closed curves $R,R_1,R_2$ are mutually disjoint, $h_R,h_{R_1},h_{R_2}$ are all commutative. We have
    \[
    \begin{aligned}
      x_1&\xmapsto{h_R'}R_1x_2R^{-1}\xmapsto{h_A}R_1x_2R_1^{-1}\\
      y_1&\xmapsto{h_R'}Ry_2^{-1}R_1^{-1}\xmapsto{h_A}R_1y_2R_1^{-1}\\
      x_2&\xmapsto{h_R'}y_1x^{-1}y_1^{-1}\xmapsto{h_A}x_1\\
      y_2&\xmapsto{h_R'}R_1y_1^{-1}\xmapsto{h_A}y_1\\
    \end{aligned}
    \]
    as desired.
  \end{proof}
\end{theorem}

\begin{figure}[hbt]
\begin{center}
  \includegraphics{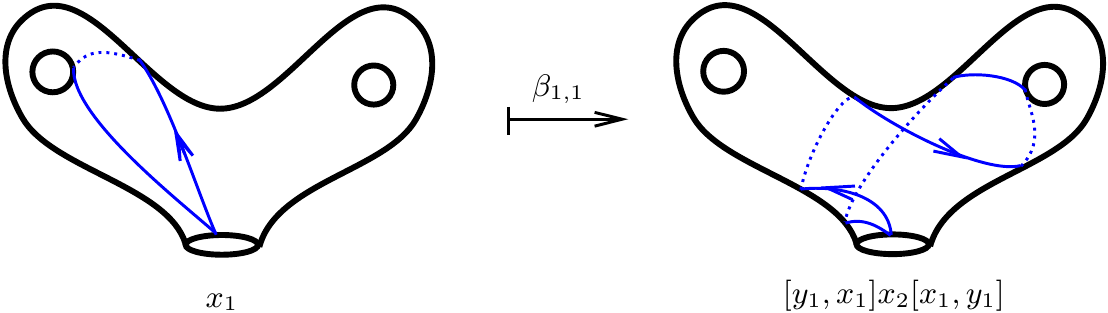}\\
  \caption{The action of $\beta_{1,1}$}\label{fig6}
\end{center}
\end{figure}

\begin{figure}[hbt]
\begin{center}
  \includegraphics{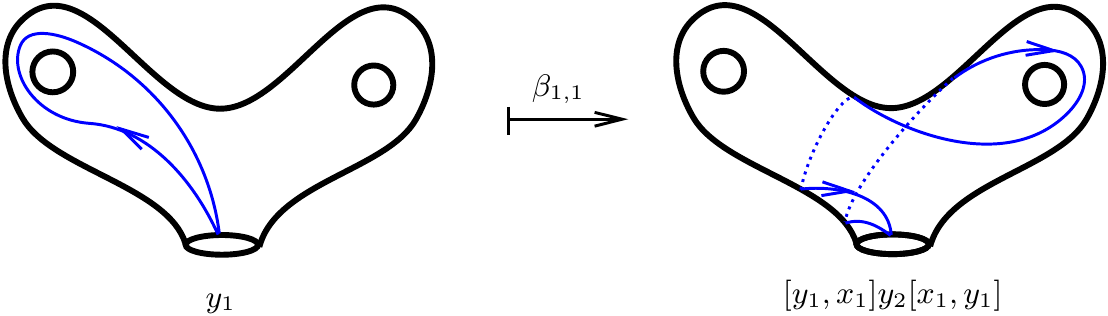}\\
  \caption{The action of $\beta_{1,1}$}\label{fig7}
\end{center}
\end{figure}

We have shown that $\beta_{1,1}$ is the composition of the reverse half Dehn twist along the boundary and the half Dehn twists on two arms. It is illustrated in Figure~\ref{fig6} and ~\ref{fig7}.

Recall that the $(r,s)$-braiding $\beta_{r,s}\in\Gamma_{r+s,1}$ equals
\[(\beta_r\beta_{r-1}\cdots\beta_1)(\beta_{r+1}\beta_r\cdots\beta_2)\cdots(\beta_{r+s-1}\beta_{r+s-2}\cdots\beta_s),\]
where $\beta_i=(a_ib_ia_i)^4(a_{i+1}b_{i+1}(a_ib_ia_i)^{-1}w_ia_ib_ia_i^2b_i)^{-3}$.  Let $R_i=[y_i,x_i]$. Then $\beta_i\in\Gamma_{r+s,1}$ acts on $\pi_1S_{r+s,1}$ as follows:
\[
\begin{array}{rlll}
  \beta_i:&x_i\mapsto R_ix_{i+1}R_i^{-1}&,&y_i\mapsto R_iy_{i+1}R_i^{-1}\\
  &x_{i+1}\mapsto x_i&,&y_{i+1}\mapsto y_i.
\end{array}
\]

\begin{figure}[hbt]
\begin{center}
  \includegraphics{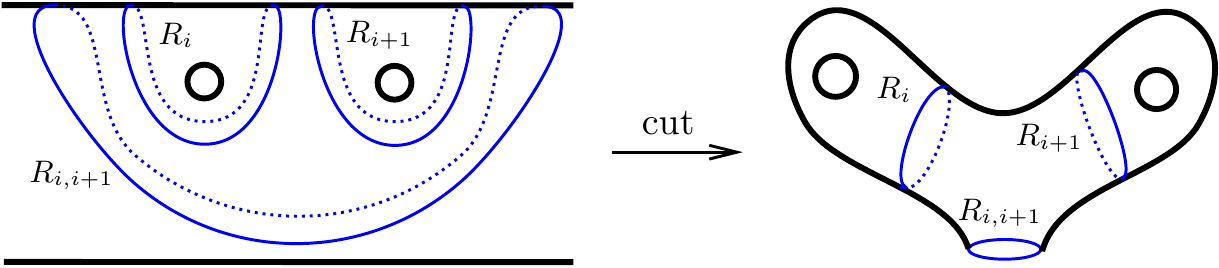}\\
  \caption{The local surface $S_{2,1}$ obtained by cutting along $R_{i,i+1}$}\label{fig8}
\end{center}
\end{figure}

Let $R_{i,j}=[y_i,x_i][y_{i+1},x_{i+1}]\cdots[y_j,x_j]$ for $j>i$, For a closed curve $R_{i,i+1}$ on $S_{r+s,1}$ (See Figure~\ref{fig8}) imagine that we cut out the surface along $R_{i,i+1}$, then we get a (small) surface $S_{2,1}$.

Then $\beta_i$ is a self-homeomorphism of $S_{r+s,1}$ which is the composition of the reverse half Dehn twist along $h_{R_{i,i+1}}'$ and the half Dehn twists on two arms, that is, we have
\[\beta_i=h_{R_{i,i+1}}'\circ h_{R_i}\circ h_{R_{i+1}}.\]

We can see, by checking the actions on the fundamental group, that these $\beta_i$'s satisfy the braid relations:
\[\beta_i\beta_{i+1}\beta_i=\beta_{i+1}\beta_i\beta_{i+1}\]
\[\beta_i\beta_j=\beta_j\beta_i\mathrm{\ for\ }|i-j|\geq 2.\]

\begin{figure}[hbt]
\begin{center}
  \includegraphics{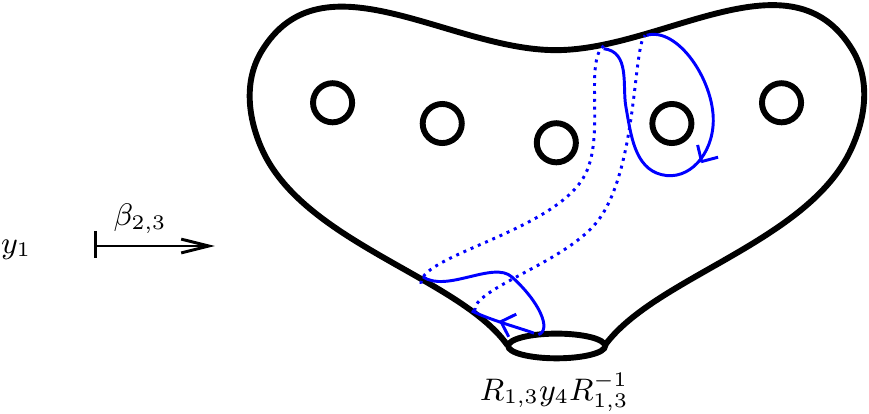}\\
  \caption{An example of action of braiding}\label{fig9}
\end{center}
\end{figure}

The $(r,s)$-braiding $\beta_{r,s}$ is a composition of a series of local braidings $\beta_{i}'s$. It may be regarded as the composition of the reverse half Dehn twist along the boundary $R$ of $S_{r+s,1}$ and a series of some local half Dehn twists in the upper part of the whole surface.  Figure~\ref{fig9} shows that how $\beta_{2,3}$ acts on $y_1$.

Let $\aut F_n$ be the automorphism of the free group $F_n$ on $n$ generators $x_1,\ldots,x_n$. Then there is an injection $j:B_n\ra\aut F_n$, called the Artin map (\cite{Artin}, \cite{Tillmann}), defined by
\[j(\sigma_i):x_i\mapsto x_ix_{i+1}x_i^{-1},\ x_{i+1}\mapsto x_i\]
where $\sigma_1,\ldots,\sigma_{n-1}$ are the standard generators of $B_n$.

There exists an obvious group homomorphism
\[
\begin{aligned}
  \phi:&B_g\ra\Gamma_{g,1}\subset\aut S_{g,1}\\
  &\sigma_i\mapsto\beta_i
\end{aligned}
\]
where $\phi(\sigma_i)=\beta_i$ acts on $\pi_1S_{g,1}=F_{\{x_1,y_1,x_2,y_2,\ldots,x_g,y_g\}}$ as
\[
\begin{array}{rll}
  x_i\mapsto R_ix_{i+1}R_i^{-1}&,&y_i\mapsto R_iy_{i+1}R_i^{-1}\\
  x_{i+1}\mapsto x_i&,&y_{i+1}\mapsto y_i,
\end{array}
\]
where $R_i=[y_i, x_i]$. The map $\phi$ is injective as the Artin map $j:B_n\ra\aut F_n$ is. (\cite{Artin}, \cite{Tillmann})

An embedding of a braid group into a mapping class group is called a geometric embedding if it takes the standard generators of braid group onto Dehn twists in the mapping class group(\cite{W2}). In \cite{W3} Wajnryb raised a question whether there can be a nongeometric embedding. We here show that map $\phi : B_g \rightarrow \Gamma_{g,1}$ is a nongeometric embedding.
\begin{lemma}
The map $\phi : B_g \rightarrow \Gamma_{g,1}$ is a nongeometric embedding.
\begin{proof}
For a standard generator $\sigma_i$ of $B_g$, $\phi(\sigma_i)=\beta_i$ is composite of three half Dehn twists along three disjoint simple closed curves in the local surface $S_{2,1}$. We may think that this self-homeomorphism is supported in the inner surface $S_{0,3}$ obtained by removing the outer dotted parts of the surface $S_{2,1}$ as given in Figure~\ref{fig_pants}.  $\beta_i$ cannot be equal to a full Dehn twist along a simple closed curve in $S_{0,3}$ because there are only three simple closed curves in $S_{0,3}$ and any Dehn twist along one of those three closed curves is not equal to $\phi(\sigma_i)$. We can see this by comparing the actions of them on the fundamental group of $S_{2,1}$.
\end{proof}
\end{lemma}

\begin{figure}[hbt]
\begin{center}
  \includegraphics{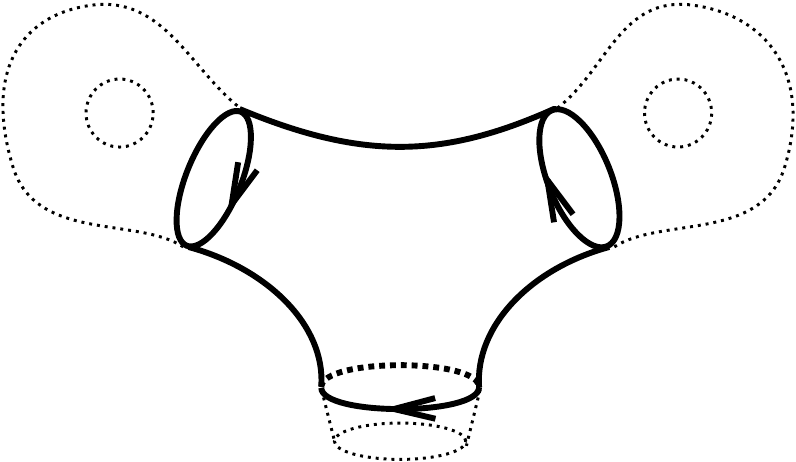}\\
  \caption{Three $180^{\circ}$ twists along the boundary of $S_{0,3}$}\label{fig10}
\end{center}
\end{figure}

From the geometric construction of $\beta_i's$ and the actions of them on the fundamental group of the surface, we have the following theorem.

\begin{theorem}
  Let $\B=\coprod_{g\geq 0}B_g$ and $\M=\coprod_{g\geq 0}\Gamma_{g,1}$. The nongeometric embedding $\phi:B_g\ra\Gamma_{g,1}, \sigma_i\mapsto\beta_i$ extends to a braided monoidal functor $\Phi:\B\ra\M$.
\end{theorem}

The map $\phi:B_g\ra\Gamma_{g,1}$ is a more natural embedding than the Harer map which is defined by
\[
h(\sigma_i)=\left\{
  \begin{array}{cl}
    b_{\frac{i+1}{2}} & \textrm{ if }i\textrm{ is odd}\\
    w_{\frac{i}{2}} & \textrm{ if }i\textrm{ is even.}\\
  \end{array}
  \right.
\]
(See Figure 1 for the definitions of $b_i$'s and $w_i$'s.), because the map $h$ does not directly give rise to a braided monoidal functor for $\B$ to $\M$. (See \cite{Song2}, Theorem 3.3)

The group completion of the classifying space of a braided monoidal category is naturally homotopy equivalent to a double loop space(\cite{Berger}, \cite{CCG}, \cite{Fie}, \cite{BFSV}). Therefore, $\Phi$ induces a map  $\bar{B\B}\ra\bar{B\M}$ of double loop spaces, where $\bar{B\B}$ and $\bar{B\M}$ denote the group completions of $B\B$ and $B\M$, respectively.

Since $\bar{B\B}\simeq\Z_+\times BB_\infty^+$ and $\bar{B\M}\simeq\Z_+\times B\Gamma_\infty^+$, we have
\begin{lemma}\label{lem_null_homotopic}
  $B\phi^+:BB_\infty^+\ra B\Gamma_\infty^+$ is a map of double loop spaces, and is null-homotopic.
\end{lemma}
It is known (\cite{S-T}, Lemma 5.3) that every double loop map $f:BB_\infty^+\ra B\Gamma_\infty^+$ is null-homotopic.

Since the plus construction does not change homology groups, Lemma~\ref{lem_null_homotopic} implies that the homology homomorphism induced by $\phi$ is trivial in the stable range.
\begin{theorem}
The homomorphism $\phi_* : H_* (B_{\infty}; \Z) \longrightarrow H_* (\Gamma_{\infty}; \Z)$ induced by $\phi$ is trivial.
\end{theorem}

By the homology stability theorem (\cite{Harer}, \cite{Ivanov}) we have that $\phi_*:H_i(B_g;\Z)\ra H_i(\Gamma_{g,1};\Z)$ is trivial for $0<i<\frac{g}{2}$.

\section*{Acknowledgements}
This research was supported by Basic Science Research Program through the National Foundation of Korea (NRF) funded by the Ministry of Education, Science and Technology (2011-0004509).

\bibliographystyle{elsarticle-num}

\begin{thebibliography}{150}

\bibitem{Artin} E. Artin, {\it Theorie der Z\"opfe}, Abh. Math. Sem. Hambur. Univ. {\bf 4} (1926), 47--72.
\bibitem{BFSV} C. Batenau, Z. Fiedorowicz, R. Schw\"{a}nzl, R. Vogt, {\it Iteratoed monoidal categories}, Adv. in Math. 176(2003), 277--349.
\bibitem{Berger} C. Berger, {\it Double loop spaces, braided monoidal categories and algebraic 3-type of space}, Contemporary Math. 227(2009), 49--66.
\bibitem{CCG} P. Carrasco, A.M. Cegarra, and A.R. Garaz\'{o}n, {\it Claasifying spaces for braided monoidal categories and lax diagrams of bicategories}, arXiv0907.0930v3(2009).

\bibitem{Fie} Z. Fiedorowicz, {\it The symmetric bar construction}, Preprint, available at http://www.math.ohio-state.edu/~fiedorowicz.




\bibitem{Harer} J. Harer, {\it Stability of the homology of the mapping class groups of orientable surfaces}, Ann. of Math. (2) {\bf 121} (1985), no. 2, 215--249.
\bibitem{Ivanov} N. V. Ivanov, {\it Stabilization of the homology of Teichm\"uller modular groups}, Algebra i Analiz {\bf 1} (1989), no. 3, 110--126; Leningrad Math. J. {\bf 1} (1990), no. 3, 675--691.
\bibitem{Song} Y. Song, {\it The braidings of mapping class groups and loop spaces}, Tohoku Math. J. {\bf 52} (2000), 309--319.
\bibitem{Song2} Y. Song, {\it The action of image of brading under the harer map}, Commun. Korea math. Soc. {\bf 21} (2006), No. 2, 337--345.

\bibitem{S-T} Y. Song and U. Tillmann, {\it Braid, mapping class groups, and categorical delooping}, Math. Ann. {\bf 339} (2007), 377--393.

\bibitem{Tillmann2} U. Tillmann, {\it On the homotopy of stable mapping class group}, Invent. Math. {\bf 130} (1997), 257--275.


\bibitem{Tillmann} U. Tillmann, {\it Artin's map in stable homology}, Bull. Lond. Math. Soc. {\bf 39} (2007), no. 6, 989--992.


\bibitem{W1} B. Wajnryb, {\it A simple representation for mapping class group of an orientable surface}, Israel Journal of Math. {\bf 45} (1983), 157--174.
\bibitem{W2}B. Wajnryb, {\it Artin groups and geometric monodromy}, Invent. Math. {\bf 138} (1999), No. 3, 563--571.
\bibitem{W3}B. Wajnryb, {\it A simple representation for mapping class group of an orientable surfaceRelations in the mapping class group}, Problems on mapping class groups and related topics, Proc. Pure Math. {\bf 74} (2006), Amer. Math. Soc., 115--120.
\end{thebibliography}

\end{document}